\documentclass[a4paper, 10pt, conference]{ieeeconf}      % Use this line for a4 paper

\IEEEoverridecommandlockouts                              % This command is only needed if 
\usepackage{pgfplots}
\usepackage{tikz}
    \usetikzlibrary{arrows.meta, positioning, calc, matrix}
    \pgfdeclarelayer{background}
    \pgfdeclarelayer{nodelayer}
    \pgfdeclarelayer{edgelayer}
    \pgfdeclarelayer{foreground}
    \pgfsetlayers{background,edgelayer,nodelayer,main,foreground}
    \tikzset{
        block/.style={rectangle, draw, thick, minimum width=2cm, minimum height=1cm, align=center},
        sum/.style={circle, draw, thick, minimum width=0.25cm, minimum height=0.25cm, align=center},        
        dot/.style={circle, draw, thick, fill=black, inner sep=1pt, align=center}
    }

\usepackage[T1]{fontenc}
\usepackage[utf8]{inputenc}
\usepackage{pgfplots}
\usepackage{grffile}
\pgfplotsset{compat=newest}
\usetikzlibrary{plotmarks}
\usetikzlibrary{arrows.meta}
\usepgfplotslibrary{patchplots}

\usepackage{layouts} 
    \printinunitsof{cm} % set the unit to cm 

\usepackage{mdframed} 
\usepackage[bigfiles]{media9}

\usepackage{amsmath,amssymb,amsfonts,mathtools}
\usepackage{bm}
\usepackage{theorem}
\usepackage{xparse}
\usepackage{graphicx}

\usepackage{subfigure}
\usepackage{pgfplots}
\pgfplotsset{compat=1.18}
\usepgfplotslibrary{fillbetween}

\newcommand{\acr}[1]{\acrshort{#1}}

\renewcommand{\vec}[1]{\mathrm{vec}\left(#1\right)}

\newcommand{\mat}[1]{\ensuremath{\begin{bmatrix} #1 \end{bmatrix}}}
\newcommand{\I}[0]{[I]} % identity matrix
\newcommand{\N}[0]{\mathbb{N}} % set of natrual numbers
\newcommand{\R}[0]{\mathbb{R}} % set of real numbers
\newcommand{\constraint}[1]{\mathcal{\MakeUppercase{#1}}} % a state, input or other constraint set
\newcommand{\norm}[1]{\left\lVert #1 \right\rVert} % norm of a vector or matrix, i.e. ||...||
\newcommand{\transpose}[0]{\top} % symbol for transpose
\newcommand{\pinv}[0]{\dagger} % pseudo inverse symbol

\newcommand{\sequence}[1]{\underline{#1}}

\newcommand\mydots{\makebox[1em][c]{.\hfil.\hfil.}}

\definecolor{mycolor1}{rgb}{0.00000,0.44700,0.74100}%
\definecolor{mycolor2}{rgb}{0.85000,0.32500,0.09800}%
\definecolor{mycolor3}{rgb}{0.92900,0.69400,0.12500}%

\usepackage{multirow}
\usepackage{graphicx}
\newtheorem{remark}{Remark}

\newtheorem{assumption}{Assumption}

\usepackage{setspace}
\usepackage{filecontents}
\newif\ifShareLatex 
\ShareLatextrue      		% if you work locally

\ifShareLatex
    \usepackage[acronym, style=alttree, shortcuts, toc=true, nomain, nonumberlist]{glossaries-extra}
    \renewcommand{\makeglossaries}{\makenoidxglossaries}
    
\else
    \usepackage[acronym,style=alttree, toc=true, shortcuts, xindy, nomain, nonumberlist]{glossaries}
    \RequirePackage[xindy]{imakeidx}
\fi

\newcommand{\Symbols}{List of Symbols}
\newcommand{\Notation}{Notation}
\newglossary{symbols}{sym}{sbl}{\Symbols}
\newglossary{notation}{not}{nt}{\Notation}
\glssetwidest{THISWIDE}     % adjust length as needed

\makeglossaries % don't remove this

\loadglsentries{gloss.aux}

\title{\LARGE \bf
Active Fault Identification and Robust Control for Unknown Bounded Faults via Volume-Based Costs
}

\author{Annalena Daniels, Johannes Teutsch, Fabian Kleindienst, Marion Leibold, Dirk Wollherr% <-this % stops a space
\thanks{All authors are with the Chair of Automatic Control Engineering (LSR), Department of Computer Engineering,
        Technical University of Munich, Theresienstr. 90, 80333 Munich, Germany
        {\tt\small \{annalena.daniels, johannes.teutsch, fabian.kleindienst,
        marion.leibold, dirk.wollherr\}@tum.de. 
        }\textit{Grammarly and ChatGPT were used to improve language and readability.}}
}

\begin{document}

\maketitle
\def\thefootnote{}\footnotetext{© 2025 IEEE. This work was accepted for IEEE SysTol 2025. Personal use of this material is permitted.  Permission from IEEE must be obtained for all other uses, in any current or future media, including reprinting/republishing this material for advertising or promotional purposes, creating new collective works, for resale or redistribution to servers or lists, or reuse of any copyrighted component of this work in other works.
}\def\thefootnote{\arabic{footnote}}
\thispagestyle{empty}
\pagestyle{empty}

%%%%%%%%%%%%%%%%%%%%%%%%%%%%%%%%%%%%%%%%%%%%%%%%%%%%%%%%%%%%%%%%%%%%%%%%%%%%%%%%

% \begin{abstract}
% This paper presents a novel approach for active fault diagnosis and parameter estimation in linear systems subject to unknown faults with bounded magnitudes in closed-loop control. The method combines set-membership identification with a cost function that encourages fault identification by taking the size of the uncertainty set of the model parameters that represent the possible faults into account. It uses ellipsoidal outer approximations of the uncertainty set to guide informative excitation when needed, while smoothly transitioning back to nominal control based on updated model parameters. Unlike many existing methods, the approach does not rely on predefined fault models and instead operates solely based on known bounds on  faulty system parameters and additive disturbances. Robust constraint satisfaction is guaranteed through a tube-based model predictive control strategy. Simulation results demonstrate faster fault identification and more accurate parameter estimation compared to passive and persistently exciting adaptive strategies.
% \end{abstract}
\begin{abstract}
This paper proposes a novel framework for active fault diagnosis and parameter estimation in linear systems operating in closed-loop, subject to unknown but bounded faults. The approach integrates set-membership identification with a cost function designed to accelerate fault identification. Informative excitation is achieved by minimizing the size of the parameter uncertainty set, which is approximated using ellipsoidal outer bounds. Combining this formulation with a scheduling parameter enables a transition back to nominal control as confidence in the model estimates increases. Unlike many existing methods, the proposed approach does not rely on predefined fault models. Instead, it only requires known bounds on parameter deviations and additive disturbances. Robust constraint satisfaction is guaranteed through a tube-based model predictive control scheme. 
Simulation results demonstrate that the method achieves faster fault detection and identification compared to passive strategies and adaptive ones based on persistent excitation constraints.
\end{abstract}

%%%%%%%%%%%%%%%%%%%%%%%%%%%%%%%%%%%%%%%%%%%%%%%%%%%%%%%%%%%%%%%%%%%%%%%%%%%%%%%
\section{INTRODUCTION}
\label{sec:introduction}

% \cmnt{
% \begin{itemize}
%     \item Story: Fault detection with finite models $\to$ often unknown (infinite) models in practice $\to$ How to do fault identification and fault tolerant control?
%     \item Related work: Bounded set of models $\to$ set membership identification $\to$ robust tube MPC $\to$ PE constraint for set reduction, and other related methods.
%     \item Contribution: fault detection and fault-tolerant control based on robust tube MPC with recursive model update and set-reduction cost function.
% \end{itemize}
% }

    % As automated systems become increasingly integrated into both industry and consumer applications, ensuring their reliability, safety, and efficiency has become a key challenge in modern engineering. 
    A well-designed and thoroughly tested feedback controller typically achieves good performance and maintains safety under nominal conditions, even when uncertainty is present. However, when malfunctions occur, a classical controller may struggle to maintain safe operation. In such cases, fault detection and identification~(\acr{FI}) (often also called fault diagnosis) and fault-tolerant control methods offer a remedy~\cite{FDReview2} as they make system changes visible and therefore safe control possible. 
    % Consequently, \acr{FI} facilitates preventive actions in terms of fault-tolerant control, i.e., adaptation of the controller to guarantee safe operation despite the fault, or (predictive) maintenance. 
    \acr{FI} methods are typically categorized into passive and active methods. Passive \acr{FI}~(\acr{PFD}) methods observe the system inputs and outputs in order to detect a fault, without influencing the system inputs. Although simple to use, \acr{PFD} often fails to detect minor faults within a sufficiently short time due to a possible lack of excitation~\cite{FDReview2}. In contrast, active \acr{FI}~(\acr{AFD}) methods manipulate the system via a modification of the control input such that even minimal faults are detected and isolated more rapidly~\cite{FDActive, xu2022minimal}. 
    
    While the active excitation of the system using \acr{AFD} methods potentially leads to faster detection of faults, these auxiliary inputs might have a negative effect on control objectives like reference tracking or satisfaction of safety constraints. Existing \acr{AFD} strategies typically assume an open-loop structure, focusing primarily on finding optimal input sequences for isolating and detecting the fault~\cite{FDReview2, scott2014input}. While being effective, this focus often neglects broader control objectives
    % such as reference tracking
    and fails to account for potential impacts that auxiliary inputs may have on overall system performance. Additionally, most \acr{AFD} approaches rely on a predefined set of a finite number of faulty system descriptions \cite{qiu2023active, Tabatabaeipour2015AFDITimeVarying}, which is impractical to obtain a priori in real-world settings. With these limitations in mind, this work aims at closing the gap between effective \acr{AFD} for unknown models and robust closed-loop control using set-membership identification~(\acr{SMI}) and a control objective that is split between reference tracking and fast FI.

    \subsubsection*{Related work}  

    Using \acr{SMI} for \acr{FI} has proven effective in previous work, e.g.,~\cite{Tabatabaeipour2012Passive} for \acr{PFD} and~\cite{Tabatabaeipour2015AFDITimeVarying, Wang2020AFDSetMembership, Wang2018OnlineAFDSetMembershipEllipsoidMovingWindows} for \acr{AFD}. These methods reduce uncertainty in system parameters, thereby enabling reliable fault detection, but they assume prior fault knowledge. Recent work has addressed unknown faults~\cite{zhang2016robust, daniels2024adaptive}, though mainly in a passive fashion. Active diagnostic input design under uncertainty has been explored by~\cite{guo2024robust, feng2024distributionally}, but despite the robustness in the detection algorithms, these methods are not integrated into feedback control schemes and lack guarantees on safety during closed-loop operation. Integrating \acr{FI} into closed-loop operation can be solved by using model predictive control (\acr{MPC}), which is well-suited as it handles both performance objectives and safety constraints~\cite{MPC1}. Robust MPC (\acr{RMPC}) and adaptive \acr{MPC} further integrate constraint satisfaction under additive and parametric uncertainty~\cite{lorenzen2019robust, teutsch2024adaptive}. Robust constraint satisfaction and parameter convergence under a persistence of excitation (PE) condition is addressed in the method by~\cite{Lu2023RobustAdatptiveMPCWithPE}, which ensures that the system is sufficiently excited over time to allow for accurate parameter estimation. Dual adaptive strategies have also gained attention. An MPC approach that finds a trade-off between performance and learning by treating uncertainty reduction as an objective was formulated by~\cite{heirung2017dual}. Parsi et al.~\cite{parsi2020robust, parsi2022explicit} embed exploration incentives into the cost or constraints for safe learning in uncertain linear systems. These approaches balance control and exploration but do not explicitly penalize the size of the uncertainty sets, and are hence not tailored for \acr{AFD}.

    \subsubsection*{Contribution} 
    In this work, we propose an \acr{AFD} and fault-tolerant \acr{RMPC} framework for linear systems subject to additive disturbances. 
    % Unlike existing \acr{AFD} approaches that rely on a discrete set of predefined fault models, our method implicitly considers a continuous bounded family of potentially faulty systems, without requiring explicit fault scenarios. Using SMI, an uncertainty set is maintained over time that captures all system models consistent with the observed data and known disturbance bounds. 
    Unlike existing \acr{AFD} approaches that rely on a discrete set of predefined fault models, our method implicitly accounts for a continuous, bounded set of potentially faulty systems, without requiring explicit fault scenarios. It only requires knowledge of the range within which the system parameters may vary. By using the SMI approach, we maintain an uncertainty set over time that captures all system models consistent with the observed data and known disturbance bounds.
    Instead of relying on PE conditions to reduce the model uncertainty and detect faults, we introduce a cost function that directly penalizes the volume of the uncertainty set, represented as an ellipsoid. This encourages informative input signals that reduce uncertainty and enable FI, while balancing the objective of tracking a desired reference trajectory. Safety constraints are enforced at all times using robust MPC~\cite{lorenzen2019robust}.

    \subsubsection*{Structure}
    The remainder of the paper is organized as follows: Section~\ref{sec:prelim} introduces the problem setup along with necessary preliminaries. Section~\ref{sec:method} presents the proposed approach, detailing the \acr{FI} scheme, the design of auxiliary inputs, and their integration into fault-tolerant control. Section~\ref{sec:simulation} provides a simulation example to illustrate the method. Finally, Section~\ref{sec:conclusion} concludes the paper.

    \subsubsection*{Notation}
    The set of integers $\{a,\,\dots,\, b\}$ is denoted as $\N_a^b$.
    With $\bm{1}_n \in\R^n$, we denote a column-vector of all ones, and the identity matrix of dimension $(n \times n)$ is written as $\I_n$. 
    % Analogously, a column-vector of $n$ zeros and an $(m \times n)$ matrix of all zeros are denoted as $\bm{0}_n$ and $\bm{0}_{m \times n}$, respectively. 
    In cases where the dimensions are clear from context, the subscript is omitted. The Kronecker product of two matrices $\bm{S}_1$, $\bm{S}_2$ is denoted by ${\bm{S}_1 \otimes \bm{S}_2}$. The weighted 2-norm of a vector~${\bm{s}\in\R^n}$ is ${\norm{\bm{s}}_{\bm{S}} = \sqrt{\bm{s}^\transpose\bm{S}\bm{s}}}$ with weight matrix~${\bm{S}\in\R^{n \times n}}$. 
    % A block-diagonal matrix is denoted with $\diag{\cdot}$. 
    With $\vec{\bm{S}}$, we denote the column-wise vectorization of the matrix $\bm{S}$. The Moore-Penrose pseudo inverse of a matrix $\bm{S}$ is $\bm{S}^\pinv$. With $(\bm{s}_{k})_{k=M}^{N}$ we denote a sequence of vectors $\bm{s}_k$ with time indices ${k\in\N_{M}^{N}}$ and $M \leq N$. Further, the subscript $l \mid k$ denotes predicted quantities $l$ time steps ahead of the time step $k$. Positive definiteness of a matrix $\bm{S}$ is denoted as $\bm{S} \succ \bm{0}$.
%____________________________________________________
\section{Problem Setup and Preliminaries}
\label{sec:prelim}
    In this section, the mathematical problem setup is introduced and the adaptive \acr{RMPC} framework proposed in~\cite{lorenzen2019robust} is outlined, which will later serve as the basis for integrating the proposed \acr{AFD} method
    %introduced in Section~\ref{sec:method}
    into a robust control scheme.
    \subsection{Problem Setup}
    We consider discrete time linear systems of the form
    \begin{equation}
        \bm{x}_{k+1} = \bm{A}\bm{x}_k + \bm{B}\bm{u}_k + \bm{w}_k \, ,
        \label{eq:system}
    \end{equation}
    where ${\bm{x}_k\in\R^n}$ is the state, ${\bm{u}_k\in\R^m}$ is the input, and ${\bm{w}_k\in\constraint{W}\subset\R^n}$ is a disturbance with compact polytopic support $\constraint{W} \coloneqq \left\{ \bm{w} \in \mathbb{R}^{n} ~\mid~ \bm{G}_w \bm{w} \le \bm{g}_w \right\}$. 
    Faults are modeled as changes in the system matrices $\bm{A}$ and $\bm{B}$ that occur at an unknown time. Rather than assuming a finite set of fault modes, all possible post-fault systems are presented as lying within a known uncertainty set of parameters. Specifically, each fault corresponds to a new but fixed pair $(\bm{A}, \bm{B})$ satisfying the following assumption.
    
    \begin{assumption}[Set of possible faults]
        For all admissible faults, the pair of matrices $(\bm{A},~\bm{B})$ from \eqref{eq:system} is controllable and satisfies $\bm{\theta} \coloneqq \mathrm{vec}\left(\mat{\bm{A} &\bm{B}}\right) \in \mathcal{AB}$, where
        \begin{equation} \label{eq:setsysmat_initial}
            \mathcal{AB} \coloneqq \left\{\bm{\theta} \in \R^{n(n+m)} ~\mid~ \bm{G}_{AB} \bm{\theta} \le \bm{g}_{AB} \right\}
        \end{equation}
        is a compact polytopic set. The nominal (fault-free) system matrices $\bm{A}_{\mathrm{nom}}, \bm{B}_{\mathrm{nom}}$ are known and satisfy $\mathrm{vec}\left(\mat{\bm{A}_{\mathrm{nom}} & \bm{B}_{\mathrm{nom}}}\right) \in \mathcal{AB}_0$, where $\mathcal{AB}_0$ denotes the initial set of consistent parameters before any fault occurs and is known, i.e., $\bm{G}_{AB,0}$ and $\bm{g}_{AB,0}$ are known, and ${\mathcal{AB}_k \subseteq \mathcal{AB}_0}$ must hold for every time step $k$. \hfill \small{$\square$} 
    \end{assumption}
    % This formulation allows for a continuous and bounded set of possible faults for system~\eqref{eq:system}, where a fault is defined as a change in the system matrices to a new value within $\mathcal{AB}$, starting at $(\bm{A},~\bm{B}) = (\bm{A}_{\mathrm{nom}}, \bm{B}_{\mathrm{nom}})$. 
    The set $\mathcal{AB}$ thus implicitly represents all potential fault realizations. Requiring controllability ensures that any such post-fault system remains stabilizable and suitable for robust control.

% \begin{assumption}[Set of possible faults]
%     For all admissible faults, the pair of matrices $(\bm{A},~\bm{B})$ from \eqref{eq:system} is controllable and satisfies $\bm{\theta} \coloneqq \mathrm{vec}\left(\mat{\bm{A} & \bm{B}}\right) \in \mathcal{AB}$, where $\mathcal{AB}$ is a compact, but not necessarily known, polytopic set representing all possible fault realizations consistent with the system dynamics~\eqref{eq:system}. 

%     An initial known outer approximation of this set, denoted $\mathcal{AB}_0$, is defined by known matrices $\bm{G}_{AB}, \bm{g}_{AB}$ as
%     \begin{equation} \label{eq:setsysmat_initial}
%         \mathcal{AB}_0 \coloneqq \left\{\bm{\theta} \in \R^{n(n+m)} ~\mid~ \bm{G}_{AB} \bm{\theta} \le \bm{g}_{AB} \right\}.
%     \end{equation}
%     This initial set contains the nominal (fault-free) system matrices, i.e., $\mathrm{vec}\left(\mat{\bm{A}_{\mathrm{nom}} & \bm{B}_{\mathrm{nom}}}\right) \in \mathcal{AB}_0$, and is used to initialize the set-membership model. At each time step $k$, the set of consistent parameters $\mathcal{AB}_k$ satisfies $\mathcal{AB}_k \subseteq \mathcal{AB}_0$. \hfill \small{$\square$}
% \end{assumption}

% This formulation captures a continuous family of possible faults, where each fault corresponds to a change in the system matrices to a new value within the (unknown) set $\mathcal{AB}$. Requiring controllability ensures that all potential fault realizations remain amenable to robust control and fault-tolerant design.

    System \eqref{eq:system} is subject to safety constraints, i.e., $\forall k \ge 0$ %for all time steps $k \ge 0$, i.e.,
    \begin{subequations}\label{eq:constraints}
    \begin{align}
    	&\bm{x}_{k} \in \constraint{X},~~ \constraint{X} = \left\{\bm{x} \in \mathbb{R}^{n}  ~\mid~ \bm{G}_x \,\bm{x} \le \bm{g}_x \right\}, \label{eq:statecons} \\
    	&\bm{u}_{k} \in  \constraint{U},~~ \constraint{U} = \left\{\bm{u} \in \mathbb{R}^{m} \,\mid~ \bm{G}_u \bm{u} \le \bm{g}_u \right\}, \label{eq:inputcons}
    	\end{align}
    \end{subequations}
    with compact polytopic sets $\constraint{X}$ and $\constraint{U}$ containing the origin. 

    % This paper aims to detect deviations from the nominal behavior $(\bm{A}{\mathrm{nom}}, \bm{B}{\mathrm{nom}})$, identify the current system matrices $(\bm{A}_k, \bm{B}_k)$ as quickly as possible, and safely pursue a control objective, all while satisfying the constraints~\eqref{eq:constraints} under uncertainty.
    The goal is to detect deviations from the nominal behavior $(\bm{A}_{\mathrm{nom}}, \bm{B}_{\mathrm{nom}})$, identify the current system matrices $(\bm{A}_k, \bm{B}_k)$ as quickly as possible, and pursue a control objective, all while satisfying the constraints~\eqref{eq:constraints} under uncertainty arising from unknown system parameters $\bm{\theta}_k \in \mathcal{AB}_k$ and disturbances $\bm{w}_k \in \mathcal{W}$ for all $k \ge 0$.

    \subsection{Adaptive Robust Model Predictive Control} \label{sec:rmpc}
    In \cite{lorenzen2019robust}, an adaptive RMPC framework with online pa\-ra\-me\-ter estimation was proposed, allowing for robust satisfaction of safety constraints despite model uncertainty and additive disturbances, which is recalled here. For now, consider system~\eqref{eq:system} with fixed but unknown system matrices $\bm{A}$, $\bm{B}$ satisfying $\mathrm{vec}\left(\mat{\bm{A} &\bm{B}}\right) \in \mathcal{AB}_k$ where $\mathcal{AB}_k$ is a time-varying compact set.

    As common in RMPC for disturbance attenuation, the control input is parameterized as 
    % \begin{equation} \label{eq:inputdecomp}
    %     \bm{u}_k = \bm{K} \bm{x}_k + \bm{v}_k,
    % \end{equation}
    $\bm{u}_k = \bm{K} \bm{x}_k + \bm{v}_k$,
    % with a stabilizing feedback gain $\bm{K}$ and the input correction term $\bm{v}_k$. Note that stabilizing feedback gains can be directly computed from a parameter set as in \eqref{eq:setsysmat_initial} by solving linear matrix inequalities, see~\cite{scherer2000linear}.
    with a stabilizing gain $\bm{K}$ and the correction term $\bm{v}_k$. Note that stabilizing gains can be directly computed from a parameter set as in \eqref{eq:setsysmat_initial} by solving linear matrix inequalities~\cite{scherer2000linear}.
    
    In order to satisfy the constraints~\eqref{eq:constraints} despite the uncertainty sets $\mathcal{AB}_k$ and $\constraint{W}$, a homothetic state tube is employed as proposed in~\cite{rakovic2012homothetic}. This state tube consists of sets $\constraint{X}_{l\mid k}$, $l\in \N_0^N$, that satisfy the constraints~\eqref{eq:constraints} over the prediction horizon $N$ for all possible choices $\bm{A}$, $\bm{B}$, $\bm{w}$, and act as an outer bound for the predicted states $\bm{x}_{l\mid k}$. That is, $\forall l\in\N_0^N$:
    \begin{subequations}
        \label{eq:htmpc_constraints}
        \begin{align}
            \bm{x}_{l\mid k} &\in \constraint{X} \hspace{8mm}\forall \bm{x}_{l\mid k} \in  \constraint{X}_{l\mid k}, \\
            \bm{K}\bm{x}_{l\mid k} + \bm{v}_{l\mid k} &\in \constraint{U} \hspace{8mm}\forall \bm{x}_{l\mid k} \in  \constraint{X}_{l\mid k}, \\
            (\bm{A} + \bm{B}\bm{K})\bm{x}_{l\mid k} + \bm{B}\bm{v}_{l\mid k} + w &\in \constraint{X}_{l+1\mid k} ~\forall \bm{x}_{l\mid k} \in \constraint{X}_{l\mid k}, \notag\\
            &\hspace{-13mm}\vec{\mat{\bm{A}&\bm{B}}} \in \mathcal{AB}_k,~ w\in\constraint{W}.
        \end{align}
    \end{subequations}
    Note that $(\bm{A},~\bm{B})$ and their bounds specified in $\mathcal{AB}_k$ are constant and not updated in the prediction at time $k$. 
    A common parameterization of the sets $\constraint{X}_{l\mid k}$ is given as $\constraint{X}_{l\mid k} \coloneqq \{\bm{z}_{l\mid k}\} \oplus \alpha_{l\mid k}\constraint{X}_{\mathrm{T}}$, where $\oplus$ denotes the Minkowski sum, $\constraint{X}_{\mathrm{T}}$ is a user-chosen set defining the shape and complexity of the tube, and $\bm{z}_{l\mid k} \in \R^n$ and $\alpha_{l\mid k} > 0$ are the center and scaling of the set, which must be determined in every time step.
    
    Thus, the optimal control problem (OCP) of the RMPC framework for time step $k$ is formulated as %delivers a control input $\bm{v}_{0\mid k}$ that is applied at time~$k$
    \begin{equation} \label{eq:ocp}
        \underset{\sequence{\bm{v}}_{N,k}}{\mathrm{minimize}} ~~ J_N(\bm{x}_k,\, \sequence{\bm{v}}_{N,k},\,\hat{\bm{A}}_k,\,\hat{\bm{B}}_k) ~~ \text{s.t.} ~~ \eqref{eq:htmpc_constraints} ~~\forall l\in\N_0^N ,
    \end{equation}
    with $\sequence{\bm{v}}_{N,k} \coloneqq \{\bm{v}_{0\mid k},\,\dots,\,\bm{v}_{N-1\mid k}\}$ and $\bm{x}_{0|k} = \bm{x}_k$. Considering $\bm{u}_{l\mid k} = \bm{K} {\bm{x}}_{l\mid k} + \bm{v}_{l\mid k}$, we define the cost function as
    \begin{equation}
        J_{N} \coloneqq \sum_{l=0}^{N-1} \norm{{\bm{x}}_{l\mid k} - \bm{x}_{k+l}^{\mathrm{ref}}}_{\bm{Q}}^2 + \norm{\bm{u}_{l\mid k} - \bm{u}_{k+l}^{\mathrm{ref}}}_{\bm{R}}^2,
        \label{eq:cost_function_tracking}
    \end{equation}
    with weighting matrices $\bm{Q},\,\bm{R} \succ \bm{0}$.
    Here, $\bm{x}_{k+l}^{\mathrm{ref}}$, $\bm{u}_{k+l}^{\mathrm{ref}}$ are references and ${\bm{x}}_{l\mid k}$ is predicted based on the input sequence $\sequence{\bm{v}}_{N,k}$ and the parameter estimate $\hat{\bm{A}}_k$, $\hat{\bm{B}}_k$, $\mathrm{vec}\left(\mat{\hat{\bm{A}}_k &\hat{\bm{B}}_k}\right) \in \mathcal{AB}_k$, using the nominal dynamics
    \begin{equation} \label{eq:system_nominal}
        \bm{x}_{l+1\mid k} = (\hat{\bm{A}}_k + \hat{\bm{B}}_k\bm{K}) \bm{x}_{l\mid k} + \hat{\bm{B}}_k \bm{v}_{l\mid k}.
    \end{equation}
    Terminal cost and constraints can be added for stability and recursive feasibility guarantees, which are out of scope.% but are not scope of this work.

    The OCP~\eqref{eq:ocp} is solved at every time step for a given measurement $\bm{x}_k$ and the control input $\bm{u}_k = \bm{K} \bm{x}_k + \bm{v}^*_{0|k}$ is applied to the system, where $\bm{v}^*_{0|k}$ is the first element of the optimal control sequence $\sequence{\bm{v}}_{N,k}^*$. 
    Details on how to 
    % choose the parameter estimate $\hat{\bm{A}}_k$, $\hat{\bm{B}}_k$ and on how to 
    cast \eqref{eq:ocp} as an efficiently solvable quadratic program are given in~\cite{lorenzen2019robust}. 
    
    To reduce the model uncertainty, SMI is used: At time step $k \ge 1$, the set $\mathcal{AB}_{k}$ is obtained by intersection of the previous set $\mathcal{AB}_{k-1}$ with set $\Delta_k$ of parameters that are consistent with the most recent measure\-ments $(\bm{u}_{k-1},\bm{x}_{k-1},\bm{x}_{k})$, i.e., $\mathcal{AB}_{k} = \mathcal{AB}_{k-1} \cap \Delta_k$ with
    \begin{equation} \label{eq:setsysmat_newdata}
        \Delta_k \coloneqq \left\{\bm{\theta} \,\Bigm|\, - \mat{\bm{x}_{k-1}\\\bm{u}_{k-1} }^{\transpose}\otimes \bm{G}_w \bm{\theta} \le \bm{g}_w - \bm{G}_w \bm{x}_{k} \right\}.
    \end{equation}
    By starting with the initial set $\mathcal{AB}_0$, the size of the set of model uncertainty monotonically decreases.
    % i.e., $\mathcal{AB}_{k+1} \subseteq \mathcal{AB}_{k}$. 
    In fact, following ideas from \cite{marafioti2014persistently}, it was shown in \cite{lorenzen2019robust} that the set $\mathcal{AB}_{k}$ converges to a singleton containing the true parameters $\bm{A}$, $\bm{B}$ when incorporating a PE constraint, with some suitable parameter $a > 0$, into the OCP \eqref{eq:ocp}, i.e., 
    \begin{equation}
        \left(\sum_{i=0}^n \bm{u}_{k-i} \bm{u}_{k-i}^{\transpose}\right) - a \bm{I}_m \succ \bm{0}.
        \label{eq:PEconstraint}
    \end{equation}

    However, the desirable properties of the SMI approach rely on the assumption that the system matrices remain constant, i.e., no fault occurs. If a fault alters the system parameters, the update law $\mathcal{AB}_{k} = \mathcal{AB}_{k-1} \cap \Delta_k$ may result in an empty set, since the collected data originate from different systems. Moreover, while condition~\eqref{eq:PEconstraint} ensures convergence, it offers no guarantees on the rate of convergence, making it less suitable for fast FI and fault-tolerant control.
    To address these limitations, a method is proposed that maintains robustness in the presence of unknown faults and supports fast FI.
    % system adaptation, constraint-aware .

% \section{Method}
% \label{sec:method}
% Our framework for the \acr{AFD} of unknown bounded faults consists of three core components. In Section~\ref{sec:faultdetection}, it is shown how \acr{SMI} is leveraged to both detect and diagnose faults in a passive fashion. Based on this, in Section~\ref{sec:auxiliary}, an active excitation strategy is proposed to accelerate the \acr{FI} process by minimizing the volume of ellipsoidal outer approximations of the uncertainty set. Unlike methods that rely on a fixed set of discrete fault models, a continuous fault space is considered here that must be progressively narrowed. This difference requires designing inputs that reduce uncertainty, rather than simply maximizing discriminability across known fault scenarios~\cite{Tabatabaeipour2015AFDITimeVarying}. Finally, Section~\ref{sec:faulttolerant} introduces the fault-tolerant control strategy based on the proposed \acr{AFD} approach.

\section{Method}
\label{sec:method}
Our framework for \acr{AFD} of unknown bounded faults consists of three core components. Section~\ref{sec:faultdetection} shows how \acr{SMI} enables passive fault detection and diagnosis. Building on this, Section~\ref{sec:auxiliary} introduces an active excitation strategy to accelerate \acr{FI} by minimizing the volume of ellipsoidal outer approximations of the uncertainty set. Unlike methods based on a fixed set of discrete fault models, we consider a continuous fault space that must be progressively narrowed. This requires designing inputs that reduce uncertainty rather than simply maximizing discriminability across known fault scenarios~\cite{Tabatabaeipour2015AFDITimeVarying}. Finally, Section~\ref{sec:faulttolerant} presents the fault-tolerant control strategy based on the proposed \acr{AFD} approach.

\subsection{Fault Detection and Diagnosis}
\label{sec:faultdetection}

SMI can be employed as a passive \acr{FD} mechanism, operating continuously in parallel with control. A fault is detected if the current model set $\mathcal{AB}_k$ no longer contains the nominal system parameters,
$\vec{\mat{\bm{A}_\mathrm{nom} & \bm{B}_\mathrm{nom}}} \notin \mathcal{AB}_k$, or if the set becomes empty, $\mathcal{AB}_k = \emptyset$, due to inconsistent measurements. Either condition indicates that the system dynamics have changed and are no longer consistent with the assumed nominal model.
Upon fault detection, if the set becomes empty, the SMI process is re-initialized using the prior uncertainty set $\mathcal{AB}_0$, as past observations are no longer representative of the current system behavior.

% After the fault is detected, the diagnosis stage begins and the fault model is estimated at each time step. 
% This can be done by using the geometric center or the Chebyshev center of the current uncertainty set $\mathcal{AB}_k$ or alternatively, a least-squares estimate based on observed input/state data.
% The least-squares estimate is obtained by minimizing the squared error
After fault detection, the diagnosis stage begins, and the fault model is estimated at each time step. This can be done using the geometric or Chebyshev center of the current uncertainty set~$\mathcal{AB}_k$, or alternatively, via a least-squares estimate (LSE) based on observed input/state data. The LSE is obtained by minimizing the squared error
\begin{equation} \label{eq:lsopt}
    \underset{\bm{A},\,\bm{B}}{\mathrm{minimize}} \ \left\| \bm{H}_{x+} - 
    \mat{\bm{A} & \bm{B}} \bm{H}_{xu} \right\|_F^2,
\end{equation}
where $\left\| \cdot \right\|_F^2$ denotes the squared Frobenius norm and 
\begin{align*}
    \bm{H}_{xu} &= \mat{\bm{x}_{k-N_{\mathrm{ls}}}\,\mydots\,\bm{x}_{k-1}\\ \bm{u}_{k-N_{\mathrm{ls}}}\,\mydots\,\bm{u}_{k-1}}, &\bm{H}_{x+} &=\mat{\bm{x}_{k-N_{\mathrm{ls}}+1}\,\mydots\,\bm{x}_{k}}
\end{align*} 
are matrices consisting of past input/state data over a user-chosen horizon $N_{\mathrm{ls}} > 0$. The closed-form solution of \eqref{eq:lsopt} is $\mat{\bm{A}_{\mathrm{est}} & \bm{B}_{\mathrm{est}}} = \bm{H}_{x+} \bm{H}_{xu}^\dag. $

As discussed in Section~\ref{sec:rmpc}, a nominal system model \eqref{eq:system_nominal} with parameters $\hat{\bm{A}}_k$, $\hat{\bm{B}}_k$ is required for the cost function of the OCP~\eqref{eq:ocp}. These models are updated as
\begin{equation}
\label{eq:main:model_adaptation}
\begin{aligned}
    (\hat{\bm{A}}_k,\, \hat{\bm{B}}_k) = 
    \begin{cases}
        (\bm{A}_{\mathrm{est}},\, \bm{B}_{\mathrm{est}}) & \text{if fault detected}, \\[0.5ex]
        (\bm{A}_{\mathrm{nom}},\, \bm{B}_{\mathrm{nom}}) & \text{otherwise,}
    \end{cases}
\end{aligned}
\end{equation}
as it is assumed that the nominal model describes the behavior best as long as this does not lead to inconsistencies. 

% While this algorithm can operate in a fully passive manner, systems with repetitive trajectories or near-equilibrium behavior may yield insufficiently informative data. In such cases, small faults may be detected late or missed entirely, and post-fault models may be inaccurate due to rank deficiency in the data matrix $\bm{H}_{xu}$ used in the least-squares estimation~\eqref{eq:lsopt}, degrading control performance. To address this, the next section discusses how to generate auxiliary inputs that excite the system and improve the \acr{FI}.

Although the algorithm can operate passively, repetitive trajectories or near-equilibrium behavior may in some cases yield uninformative data. This can delay or prevent fault detection and lead to inaccurate post-fault models due to rank deficiency in the data matrix $\bm{H}_{xu}$ used in the LSE~\eqref{eq:lsopt}, degrading control performance. To address this, the next section introduces auxiliary inputs to excite the system and improve \acr{FI}.

\subsection{Design of Safe Auxiliary Inputs}
\label{sec:auxiliary}
Encouraging informative inputs that reduce model uncertainty can be intuitively achieved by incorporating the volume of the model uncertainty set into the MPC cost function. However, computing the volume of a parametric uncertainty set, typically a polytope, requires solving a separate optimization problem. Embedding this into the MPC cost results in a nested (bilevel) optimization structure, which is generally intractable for online use due to its nonconvexity and high computational cost.

To overcome this challenge, an ellipsoidal outer approximation of the uncertainty set $\bm{\theta} = \mathrm{vec}([\bm{A} \ \bm{B}])$ is used. Specifically, the model uncertainty is represented by% an ellipsoid
\begin{equation}
    \mathcal{T} = \left\{ \bm{\theta} \in \mathbb{R}^{n(n + m)} \ \middle| \ (\bm{\theta} - \bm{c})^{\transpose} \bm{C}^{-1} (\bm{\theta} - \bm{c}) \leq 1 \right\},
    \label{eq:ellispoid}
\end{equation}
where \( \bm{c} \in \mathbb{R}^{n(n + m)} \) is the ellipsoid center, and \( \bm{C} \succ 0 \) is the symmetric positive definite shape matrix. The volume \( V \) of this ellipsoid is proportional to \( \det(\bm{C})^{1/2} \)~\cite[Chap. 8.4]{boyd2004convex}, which is equivalent to \( V \propto \det(\bm{C}^{-1})^{-1/2} \). We incorporate a convex cost term based on this volume into the MPC objective by penalizing \( -\log\det(\bm{C}^{-1}) \), thereby discouraging large uncertainty. In order to find $\bm{C}^{-1}$, the ellipsoidal approximation of $\bm{\theta}$ needs to be constructed. The system equations are rewritten in a stacked form over a batch of \( k+N \) data points
\begin{equation} \label{eq:disturbanceeq}
    \sequence{\bm{y}}_{k+N} = \bm{Z}_{k+N} \bm{\theta} + \sequence{\bm{w}}_{k+N},
\end{equation}
with
%\begin{subequations}
\begin{align}
    \sequence{\bm{y}}_{k+N}  &= \mat{\bm{x}_1^{\transpose},\,\mydots,\,\bm{x}_k^{\transpose},\,\bm{x}_{1\mid k}^{\transpose},\,\mydots,\,\bm{x}_{N\mid k}^{\transpose}}^{\transpose}, \notag \\
    \sequence{\bm{w}}_{k+N} &= \mat{\bm{w}_0^{\transpose},\,\mydots,\,\bm{w}_{k-1}^{\transpose},\,\bm{w}_{k}^{\transpose},\,\mydots,\,\bm{w}_{k+N-1}^{\transpose}}^{\transpose}, \label{eq:stackedparams}\\
    \bm{Z}_{k+N} &= \mat{\bm{x}_0,\,\mydots,\,\bm{x}_{k-1},\,\bm{x}_{k},\,\mydots,\,\bm{x}_{N-1\mid k} \\ \bm{u}_0,\,\mydots,\,\bm{u}_{k-1},\,\bm{u}_{0\mid k},\,\mydots,\,\bm{u}_{N-1\mid k}}^{\transpose} \otimes I_n. \notag
\end{align}

%\end{subequations}
The vector $\sequence{\bm{y}}_{k+N}$ and the matrix $\bm{Z}_{k+N}$ consist of past inputs and states up to time step $k$, and of predicted inputs $\bm{u}_{l\mid k} = \bm{K} {\bm{x}}_{l\mid k} + \bm{v}_{l\mid k}$ and states using the nominal prediction model \eqref{eq:system_nominal}. The vector $\sequence{\bm{w}}_{k+N}$ consists of the corresponding disturbances, which are unknown but satisfy the bounds $\bm{w}_i \in \constraint{W} ~\forall i \in \N_0^{k+N}$. As the polytopic disturbance bound is known, an ellipsoid outer approximation $\constraint{W}_{k+N}$ of the disturbance set can be derived that contains $\sequence{\bm{w}}_{k+N}$, i.e.,
\begin{equation}
    \constraint{W}_{k+N} = \left\{\sequence{\bm{w}} \ \middle|\ (\sequence{\bm{w}} - \bm{c}_w)^{\transpose} \bm{C}_w^{-1} (\sequence{\bm{w}} - \bm{c}_w) \leq 1 \right\},
    \label{eq:disturbanceellips}
\end{equation}
where \( \bm{c}_w \in \mathbb{R}^{n(k+N)}  \) is the center and \( \bm{C}_w \succ 0 \) the shape matrix of the ellipsoid. 
Computing the smallest outer approximation \eqref{eq:disturbanceellips}, ensuring that all disturbances $\bm{w}_i \in \constraint{W} ~\forall i \in \N_0^{k+N}$ lie within the set, leads to a convex optimization problem that can be solved offline. Efficient algorithms such as Khachiyan’s ellipsoid method or semidefinite programming formulations can be used to compute a minimum-volume enclosing ellipsoid; see, e.g.,~\cite{kumar2005minimum}.
% \begin{remark}
%     One could consider modeling the disturbance set as an ellipsoid from the beginning, thereby avoiding the need for ellipsoidal outer approximations during volume computation. A robust MPC framework based on ellipsoids was presented in \cite{parsi2025scalable}; however, it does not incorporate parameter estimation or adaptation. Extending such an ellipsoidal approach to include parameter adaptation would mean representing the parameter uncertainty set as an ellipsoid as well. While this would allow the use of convex tools and maintain convexity under intersection, the resulting intersection of multiple ellipsoids is not itself an ellipsoid. Therefore, additional outer approximations would be needed after each update step to retain an ellipsoidal structure. Studying the implications of this trade-off and designing a fully ellipsoidal identification and control framework is a research topic in its own right and lies beyond the scope of this work.
% \end{remark}

By solving \eqref{eq:disturbanceeq} for $\sequence{\bm{w}}_{k+N}$ and leveraging its bounds \eqref{eq:disturbanceellips}, one obtains that any parameter $\bm{\theta}$ consistent with the past data and future predictions must satisfy
\begin{equation}
    (\sequence{\bm{y}}_{k+N} - \bm{Z}_{k+N} \bm{\theta} - \bm{c}_w)^{\transpose} \bm{C}_{w}^{-1} (\sequence{\bm{y}}_{k+N} - \bm{Z}_{k+N} \bm{\theta} - \bm{c}_w) \leq 1.
\end{equation}
Comparing this to \eqref{eq:ellispoid}, we find that
$ \displaystyle{\bm{C}^{-1} = \frac{\bm{Z}_{k+N}^{\transpose}\bm{C}_{w}^{-1}\bm{Z}_{k+N}}{1-\alpha}},
$
with 
$
    \alpha =  (\sequence{\bm{y}}_{k+N} -\bm{c}_w)^{\transpose}(\bm{C}_{w}^{-1} -  \tilde{\bm{Z}}^{\transpose} \bm{C}_{w}^{-1}\tilde{\bm{Z}})(\sequence{\bm{y}}_{k+N} -\bm{c}_w)
$ and $\tilde{\bm{Z}} = \bm{Z}_{k+N}\bm{Z}_{k+N}^\dag$.
% Using this approximation, we define a volume-based cost
% \begin{equation} \label{eq:setreductioncost}
%     J_{\text{vol}} = -\log\det(\bm{C}^{-1})
% \end{equation}
% as a surrogate for the size of the uncertainty set. As can be seen in \eqref{eq:stackedparams}, the data up to the current time step~$k$ contains the information of the current model uncertainty set~$\mathcal{AB}_k$ (cf. SMI approach in Section~\ref{sec:rmpc}), which is constant at time step~$k$. Thus, the cost~\eqref{eq:setreductioncost} is optimized with respect to the predicted inputs and states in~\eqref{eq:stackedparams} over the prediction horizon.
\begin{remark}
    Modeling all sets as ellipsoids would avoid the need for outer approximations during volume computation. While ellipsoidal RMPC frameworks exist~\cite{parsi2025scalable}, they typically do not include parameter adaptation. In adaptive settings, intersecting ellipsoidal uncertainty sets yields convex but non-ellipsoidal sets, which are challenging to handle directly.
    \vspace{-2mm}
\end{remark}
\vspace{-2mm}

Using this approximation, we define a volume-based cost
\begin{equation} \label{eq:setreductioncost}
J_{\text{vol}} = -\log\det(\bm{C}^{-1})
\end{equation}
as a surrogate for the size of the uncertainty set. As shown in \eqref{eq:stackedparams}, the current uncertainty set~$\mathcal{AB}_k$ is inferred from data up to time step~$k$ and therefore remains constant during prediction. The predicted inputs and states over the horizon are based on the current parameter estimates~$\hat{\bm{A}}_k$ and~$\hat{\bm{B}}_k$, and cost~\eqref{eq:setreductioncost} is optimized with respect to these predictions.

In order to enable fault identification within a controlled system, the volume-penalizing cost is incorporated into the standard RMPC formulation~\eqref{eq:ocp}. However, to avoid unnecessary excitation when identification is not required or desired, we introduce a scheduling parameter $\beta$ that regulates the influence of the volume-based cost.
This parameter allows the controller to smoothly transition between standard RMPC behavior and active uncertainty reduction, depending on the current uncertainty set volume $V$ which can be computed by directly using the ellipsoidal volumes or any other method that determines the set volume. The scaled weight $\beta(V) \in [0, 1]$ is then computed as
\begin{equation}
    \beta(V) = \log\left( \frac{\min\left( \max(V, V_{\min}), V_{\max} \right)}{V_{\min}} \right)\log\left( \frac{V_{\max}}{V_{\min}} \right)^{-1},
\end{equation}
with positive scalars $V_{\min}$ and $V_{\max}$, $0 < V_{\min} < V_{\max}$.
% Fig.~\ref{fig:beta_scaling} illustrates the scaling function $\beta(V)$, which increases from 0 to 1 as the uncertainty set volume grows. When the uncertainty is low ($V \approx V_{\min}$), the controller prioritizes reference tracking or stabilization. As the volume increases, the controller gradually shifts toward generating informative inputs to reduce model uncertainty. The tuning parameters $V_{\min}$ and $V_{\max}$ must be selected by the user based on the application. The results in~\cite{xu2022minimal} on minimal detectable and isolable faults can help inform the lower bound $V_{\min}$. If it is chosen too large, faults below a certain magnitude cannot be reliably detected as the excitation term is prematurely disabled by causing $\beta(V) \approx 0$. 
% % Thus, $V_{\min}$ should reflect the smallest detectable fault. 
% The upper bound $V_{\max}$ can be chosen based on the initial size of the uncertainty set. If falling back to pure \acr{PFD} is undesired, $\beta(V) > 0$ can be chosen.

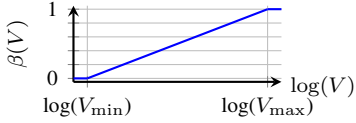
\begin{figure}[h]
    \centering
    \vspace*{-2mm}
    \begin{tikzpicture}
        \begin{axis}[
            width=0.5\linewidth,
            height=0.3\linewidth,
            axis lines=left,
            xlabel={\footnotesize $\log(V)$},
            xlabel style={at={(axis description cs:1.2,0.2)}, anchor=north},
            ylabel={\footnotesize $\beta(V)$},
            ymin=-0.05, ymax=1.05,
            xmin=-14, xmax=1,
            xtick={-13, 0},
            xticklabels={\footnotesize log($V_\mathrm{min}$), \footnotesize log($V_\mathrm{max}$)},
            ytick={0, 0.2, 0.4, 0.6, 0.8, 1},
            yticklabels={\footnotesize 0, , , , , \footnotesize 1},
            ymajorgrids=true,
            xmajorgrids=true,
            thick,
            domain=-14:1,
            samples=300
        ]
            \addplot[blue, thick, domain=-14:-13] {0};
            \addplot[blue, thick, domain=-13:0] {(x + 13) / 13};
            \addplot[blue, thick, domain=0:1] {1};
        \end{axis}
    \end{tikzpicture}
    \vspace*{-2mm}
    \caption{Logarithmic scaling of $\beta(V)$ between $V_\mathrm{min}$ and $V_\mathrm{max}$.}
    \vspace{-0.2cm}
    \label{fig:beta_scaling}
\end{figure}

Fig.~\ref{fig:beta_scaling} shows the scaling function $\beta(V)$, which increases from 0 to 1 as the uncertainty set volume grows. When $V \approx V_{\min}$, the controller prioritizes tracking or stabilization; for larger $V$, it shifts toward generating informative inputs to reduce model uncertainty. The tuning parameters $V_{\min}$ and $V_{\max}$ should reflect the application's requirements. Results from~\cite{xu2022minimal} on minimal detectable faults can inform $V_{\min}$: if chosen too large, small faults may go undetected as $\beta(V) \approx 0$ suppresses excitation. The initial uncertainty volume which can be computed offline can be used for determining $V_{\max}$. If pure \ac{PFD} is undesired, one may enforce $\beta(V) > 0$.

Consequently, the cost function $J_N$ from \eqref{eq:ocp} is adapted to
\begin{align}
    % J_N\left(\bm{x}_{k},\, \sequence{\bm{v}}_k,\, \hat{\bm{A}},\,\hat{\bm{B}},\,V \right)
    % &= J\left(\bm{x}_{k},\, \sequence{\bm{v}}_k,\, \hat{\bm{A}},\,\hat{\bm{B}} \right) \notag\\
    % &~~+ \beta(V) \cdot J_{\text{vol}}\left(\bm{x}_{k},\, \sequence{\bm{v}}_k,\, \hat{\bm{A}},\,\hat{\bm{B}} \right)
    J_N(\bm{x}_{k},\, \sequence{\bm{v}}_k,\, \hat{\bm{A}},\,\hat{\bm{B}},\,V ) &= J_{\mathrm{ctrl}}  + \beta(V) \cdot J_{\text{vol}}
    \label{eq:costtotal}
\end{align}
where $J_{\mathrm{ctrl}}$ corresponds to the nominal tracking or stabilization cost ($J_N$ in \eqref{eq:cost_function_tracking}), and $J_{\text{vol}}$ from \eqref{eq:setreductioncost} penalizes the predicted uncertainty set volume based on past observed data. 
The constraints of the homothetic tube RMPC remain unchanged. As a result, any excitation introduced for identification remains within the bounds of the tightened constraint tubes, ensuring stability, safety, and recursive feasibility. For a detailed discussion of these properties, we refer to~\cite{lorenzen2019robust}.

\subsection{Fault-tolerant Control}
\label{sec:faulttolerant}
Since the control objective is already maintained throughout the fault detection and diagnosis part, fault-tolerant control follows naturally: the algorithm simply continues operating with the updated model. Once the uncertainty set volume reaches the target threshold $V = V_\text{min}$ and an accurate model estimate $(\hat{\bm{A}}, \hat{\bm{B}}) \in \mathcal{AB}_0$ is identified, the system can be safely controlled using the adaptive RMPC framework with $\beta = 0$. If, however, the fault diagnosis mechanism finds that new measurements are no longer consistent with $\mathcal{AB}_0$, i.e., $(\hat{\bm{A}}, \hat{\bm{B}}) \notin \mathcal{AB}_0$, this indicates a more severe and unmodeled fault. In such cases, safety guarantees no longer apply, and the system should be shut down immediately.

%%%%%%%%%%%%%%%%%%%%%% Simulation and Results %%%%%%%%%%%%%%%%%%%%%%%%%%%%%%
\section{Simulation and Results}

\label{sec:simulation}
\subsection{Simulation Setup}

A discrete-time linear system with nominal dynamics
\begin{equation}
    A = \mat{A_{11} & A_{12} \\ A_{21} & A_{22} } = \mat{1 & 1 \\ 0 & 1 }, \quad
    B = \mat{ B_1 \\ B_2} =\mat{ 0.1 \\ 1},
\end{equation}
is considered. Model uncertainty is captured through a hyperrectangular initial model set $\mathcal{AB}_0$, with element-wise bounds
$ A_{11} \in [-0.8, 1.3], A_{12} \in [-0.8, 1.2], \quad A_{21} \in [-0.2, 0.2], \quad A_{22} \in [-0.8, 1.2],
    B_1    \in [-0.1, 0.2], \quad B_2    \in [0.8, 1.1]$.
State and input constraints are defined by
$    \constraint{X} = \{ \bm{x} \in \mathbb{R}^2 \mid \|\bm{x}\|_\infty \leq 5 \}, \quad
    \constraint{U} = \{ u \in \mathbb{R} \mid |u| \leq 5 \}.$
% Two different scenarios are considered in one simulation where the simulation is started with the first scenario and switches to the second one after 10 time steps. The first one is a very small fault, which is usually hard to detect and a second one with a more severe fault. The system is subject to a structural change defined by
% % \begin{equation}
% %     A_\text{fault1} = \begin{bmatrix} 1 & 0.99 \\ 0 & 1 \end{bmatrix}, \quad
% %     B_\text{fault1} = B,
% % \end{equation}
% % and
% % \begin{equation}
% %     A_\text{fault2} = \begin{bmatrix} 0.8 & 0.8 \\ 0.1 & 0.9 \end{bmatrix}, \quad
% %     B_\text{fault2} =  \begin{bmatrix} 0.15 \\ 0.95 \end{bmatrix},
% % \end{equation}
% \begin{equation}
%     A_\text{fault1} = \begin{bmatrix} 1 & 0.99 \\ 0 & 1 \end{bmatrix}, \quad
%     A_\text{fault2} = \begin{bmatrix} 0.8 & 0.8 \\ 0.1 & 0.9 \end{bmatrix},
% \end{equation}
% and $B_\text{fault1} =  B$, $B_\text{fault2} =  [0.15,\, 0.95]^{\transpose}$
% respectively. 
Two fault scenarios are considered within a single simulation: the system initially operates under a mild fault (Fault~1), which is typically difficult to detect, for 10 time steps before switching to a more severe fault (Fault~2). 
% The first scenario involves a small structural change, typically difficult to detect, while the second introduces a bigger deviation. 
The system dynamics under both faults are defined by
\begin{equation}
    A_\text{fault1} = \begin{bmatrix} 1 & 0.99 \\ 0 & 1 \end{bmatrix}, \quad
    A_\text{fault2} = \begin{bmatrix} 0.8 & 0.8 \\ 0.1 & 0.9 \end{bmatrix},
\end{equation}
with $B_\text{fault1} = B$ and $B_\text{fault2} = [0.15,\, 0.95]^{\transpose}$, respectively. The process noise is modeled as an unknown but bounded disturbance $\bm{w} \in \constraint{W} = \left\{ \bm{w} \in \mathbb{R}^{2} \ \middle| \ \| \bm{w} \|_\infty \leq 0.01 \right\}$. Disturbance samples are uniformly distributed within this set. Control is performed using the homothetic tube-based RMPC controller with a prediction horizon of $N = 3$ and $\bm{x}^\text{ref} = [0 ,\, 0]^{\transpose}$. The cost function~\eqref{eq:costtotal} uses $Q = 0.1 [I]_2$ and $R = 0.01$. %, and $Q_T = [I]_2$. 
% The terminal constraint is chosen as the robust control invariant set computed offline, details in \cite{lorenzen2019robust}. 
The initial state is $\bm{x}_0 = [0.01,\, -0.01]^{\transpose}$. The control gain is a stabilizing controller for $\mathcal{AB}_0$ and is computed to be $K = [-0.0205,\,-0.1916]$. The initial volume of $\mathcal{AB}_0$ is 0.3024, and the $\beta$ scaling parameters are set to $V_\text{min} = 10^{-13}$ and $V_\text{max}=1$. The parameters of $(\bm{A}_\text{est},~\bm{B}_\text{est})$ are the LSEs~\eqref{eq:lsopt} once matrix $\bm{H}_{xu}$ has full row-rank. Otherwise, the geometric center of $\mathcal{AB}_k$ is used.  In cases where an empty set $\mathcal{AB}_k$ is detected, $\bm{H}_{xu}$ is reset to an empty matrix.

% To compare our proposed method, three simulation scenarios are considered. In the first one, a nominal RMPC controller without excitation or volume penalization and system adaptation is used. The fault diagnosis operates in a fully passive fashion, see Section~\ref{sec:faultdetection}. The second one is the adaptive approach from~\cite{lorenzen2019robust} and includes a PE constraint as described in~\eqref{eq:PEconstraint} with $a=3$, and the third one corresponds to the proposed approach, which uses the mixed cost function~\eqref{eq:costtotal} and scaling through $\beta(V)$. Each scenario is simulated in MATLAB for 25 time steps and repeated across 120 Monte Carlo runs.

To evaluate the proposed method, three simulation scenarios are considered. The first uses a nominal RMPC controller without excitation, volume penalization, or system adaptation; fault diagnosis operates passively (see Section~\ref{sec:faultdetection}). The second follows the adaptive approach from~\cite{lorenzen2019robust}, using the PE constraint in~\eqref{eq:PEconstraint} with $a = 3$. The third implements the proposed method with the mixed cost~\eqref{eq:costtotal} and scaling via $\beta(V)$. Each scenario is simulated in MATLAB for 25 time steps and repeated over 120 Monte Carlo runs.

\subsection{Results and Discussion}

Fig.~\ref{fig:plots}a) shows the input trajectories of 5 runs for the three control strategies. The passive controller generates small inputs when the fault is mild (first 10 time steps) but fails to perfectly track the reference due to the unmodeled fault, which is reflected in the input reactions and fluctuations. This effect is even stronger for the severe fault, where tracking performance deteriorates further because the controller does not adapt to the changed system. The adaptive approach improves tracking by updating the internal model to match the real system, which results in smaller input variations even when the PE constraint is active. In contrast, the proposed \acr{AFD} method produces the largest initial excitation after the fault onsets. This behavior results from the volume-based cost term, which prioritizes informative inputs when model uncertainty is high. As the uncertainty set size decreases, the inputs, and thus the excitation, get smaller, ultimately outperforming the tracking performance of the other strategies. 
% The PE-constrained strategy, in particular, highlights the trade-off between informative excitation (high PE gain) and tracking performance (favored by low PE gains), especially since the constraint remains active and continuously enforces excitation.
\def\plotheight{0.15\textwidth}
\def\plotheighti{0.148\textwidth}

\def\plotwidth{0.31\textwidth}

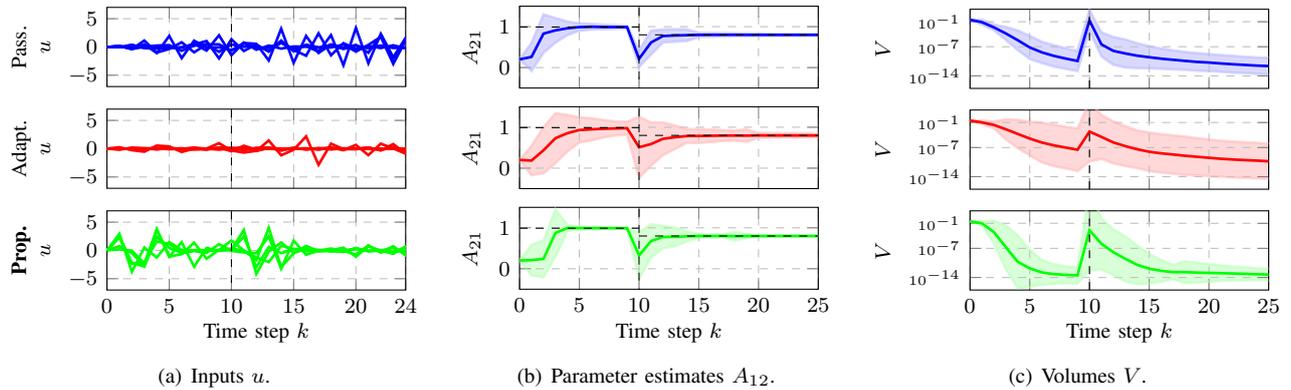
\begin{figure*}[!h]
    \vspace{0.5mm}
    \centering
    \subfigure[Inputs $u$.]{
    \begin{tikzpicture}
        % Nominal inputs:
        \begin{axis}[
            name = input_nom,
            width=\plotwidth, height=\plotheighti,
            xlabel style={at={(axis description cs:0.5,-0.3)},anchor=north},
            ylabel={Pass.\\ $u$},
            ylabel near ticks,
            xmin=0, xmax=24,
            ymin=-7, ymax=7,
            xtick={0,5,10,15,20,24},
            ytick={-5,0,5},
            xticklabels={},
            xmajorgrids=true,
            ymajorgrids=true,
            grid style=dashed,
            fill between/on layer=axis background,
            label style={font=\footnotesize, align=center},
            tick label style={font=\footnotesize},
            title style = {font=\footnotesize}
        ]

            %Nominal: (confidence bounds)
            % \addplot[name path=u_nom_ucb, line width=1pt, mark=none, blue!30, opacity=0.5, forget plot] table [x=t, y=u_nom_ucb, col sep=comma] {figs/data_means.csv};
            % \addplot[name path=u_nom_lcb, line width=1pt, mark=none, blue!30, opacity=0.5, forget plot] table [x=t, y=u_nom_lcb, col sep=comma] {figs/data_means.csv};
            % \addplot[blue!30, forget plot, opacity=0.5] fill between [of=u_nom_lcb and u_nom_ucb];
            % \addplot[line width=1pt, mark=none, blue] table [x=t, y=u_nom_mean, col sep=comma] {figs/data_means.csv};            

            %Nominal: (individual lines)
            \foreach \i in {1, 5, 8, 10, 91}
                \addplot[line width=1pt, mark=none, blue] table [x=t, y=u_nom_\i, col sep=comma] {figs/data_all2.csv};
            \addplot[black, thin, dashed] coordinates {(10, -17) (10, 7)};

        \end{axis}

        %PE inputs:
        \begin{axis}[
            name = input_pe,
            at=(input_nom.below south), anchor=above north,
            width=\plotwidth, height=\plotheighti,
            xlabel style={at={(axis description cs:0.5,-0.3)},anchor=north},
            ylabel={Adapt.\\ $u$},
            ylabel near ticks,
            xmin=0, xmax=24,
           ymin=-7, ymax=7,
            xtick={0,5,10,15,20,24},
            xticklabels={},
            ytick={-5,0,5},
            xmajorgrids=true,
            ymajorgrids=true,
            grid style=dashed,
            fill between/on layer=axis background,
            label style={font=\footnotesize, align=center},
            tick label style={font=\footnotesize} 
        ]
            
            %PE constraint: (confidence bounds)
            % \addplot[name path=u_pe_ucb, line width=1pt, mark=none, red!30, opacity=0.5, forget plot] table [x=t, y=u_pe_ucb, col sep=comma] {figs/data_means.csv};
            % \addplot[name path=u_pe_lcb, line width=1pt, mark=none, red!30, opacity=0.5, forget plot] table [x=t, y=u_pe_lcb, col sep=comma] {figs/data_means.csv};
            % \addplot[red!30, forget plot, opacity=0.5] fill between [of=u_pe_lcb and u_pe_ucb];
            % \addplot[line width=1pt, mark=none, red] table [x=t, y=u_pe_mean, col sep=comma] {figs/data_means.csv};

            %PE constraint: (individual lines)
            \foreach \i in {1, 5, 8, 10, 91}
                \addplot[line width=1pt, mark=none, red] table [x=t, y=u_pe_\i, col sep=comma] {figs/data_all2.csv};
            \addplot[black, thin, dashed] coordinates {(10, -17) (10, 7)};

        \end{axis}

        %Proposed inputs:
        \begin{axis}[
            name = input_prop,
            at=(input_pe.below south), anchor=above north,
            width=\plotwidth, height=\plotheighti,
            xlabel style={at={(axis description cs:0.5,-0.3)},anchor=north},
            ylabel={\textbf{Prop.}\\ $u$},
            xlabel={Time step $k$},
            ylabel near ticks,
            xmin=0, xmax=24,
           ymin=-7, ymax=7,
            xtick={0,5,10,15,20,24},
            ytick={-5,0,5},
            xmajorgrids=true,
            ymajorgrids=true,
            grid style=dashed,
            fill between/on layer=axis background,
            label style={font=\footnotesize, align=center},
            tick label style={font=\footnotesize} 
        ]      
            
            %Proposed: (confidence bounds)
            % \addplot[name path=u_prop_ucb, line width=1pt, mark=none, green!30, opacity=0.5, forget plot] table [x=t, y=u_prop_ucb, col sep=comma] {figs/data_means.csv};
            % \addplot[name path=u_prop_lcb, line width=1pt, mark=none, green!30, opacity=0.5, forget plot] table [x=t, y=u_prop_lcb, col sep=comma] {figs/data_means.csv};
            % \addplot[green!30, forget plot, opacity=0.5] fill between [of=u_prop_lcb and u_prop_ucb];
            % \addplot[line width=1pt, mark=none, green] table [x=t, y=u_prop_mean, col sep=comma] {figs/data_means.csv};

            %Proposed: (individual lines)
            \foreach \i in {1, 5, 8, 10, 91}
                \addplot[line width=1pt, mark=none, green] table [x=t, y=u_prop_\i, col sep=comma] {figs/data_all2.csv};
            \addplot[black, thin, dashed] coordinates {(10, -17) (10, 7)};

        \end{axis}
        
    \end{tikzpicture}
    }
    %%%%%%%%%%%%%%%%%%%%%%%%%%%%%%%%%%%%%%%%%%%%%%%%%%%%%%%%%%%%%%%%%%%%%%%%%%%
    \subfigure[Parameter estimates $A_{12}$.]{
    \begin{tikzpicture}
        % Nominal params:
        \begin{axis}[
            name = params_nom,
            width=\plotwidth, height=\plotheight,
            xlabel style={at={(axis description cs:0.5,-0.3)},anchor=north},
            ylabel={$A_{21}$},
            ylabel near ticks,
            xmin=0, xmax=25,
            ymin=-0.5, ymax=1.5,
            xtick={0,5,10,15,20,25},
            ytick={0,1},
            xticklabels={},
            xmajorgrids=true,
            ymajorgrids=true,
            grid style=dashed,
            fill between/on layer=axis background,
            label style={font=\footnotesize, align=center},
            tick label style={font=\footnotesize},
            title style = {font=\footnotesize}
        ]

            %Nominal: (confidence bounds)
            \addplot[name path=param_nom_ucb, line width=1pt, mark=none, blue!30, opacity=0.5, forget plot] table [x=t, y=param_nom_ucb, col sep=comma] {figs/data_means2.csv};
            \addplot[name path=param_nom_lcb, line width=1pt, mark=none, blue!30, opacity=0.5, forget plot] table [x=t, y=param_nom_lcb, col sep=comma] {figs/data_means2.csv};
            \addplot[blue!30, forget plot, opacity=0.5] fill between [of=param_nom_lcb and param_nom_ucb];
            \addplot[line width=1pt, mark=none, blue] table [x=t, y=param_nom_mean, col sep=comma] {figs/data_means2.csv};
            % Thin horizontal line at y = 0.99
            \addplot[black, thin, dashed] coordinates {(0, 0.99) (10, 0.99)};
            
            % Thin horizontal line at y = 0.8
            \addplot[black, thin, dashed] coordinates {(10, 0.8) (25, 0.8)};
            \addplot[black, thin, dashed] coordinates {(10, -17) (10, 2)};

            %Nominal: (individual lines)
            % \foreach \i in {1,...,10}
            %     \addplot[line width=1pt, mark=none, blue] table [x=t, y=param_nom_\i, col sep=comma] {figs/data_all.csv};
            
        \end{axis}

        %PE params:
        \begin{axis}[
            name = params_pe,
            at=(params_nom.below south), anchor=above north,
            width=\plotwidth, height=\plotheight,
            xlabel style={at={(axis description cs:0.5,-0.3)},anchor=north},
            ylabel={$A_{21}$},
            ylabel near ticks,
            xmin=0, xmax=25,
            ymin=-0.5, ymax=1.5,
            xtick={0,5,10,15,20,25},
            xticklabels={},
            ytick={0,1},
            xmajorgrids=true,
            ymajorgrids=true,
            grid style=dashed,
            fill between/on layer=axis background,
            label style={font=\footnotesize, align=center},
            tick label style={font=\footnotesize} 
        ]
            
            %PE constraint: (confidence bounds)
            \addplot[name path=param_pe_ucb, line width=1pt, mark=none, red!30, opacity=0.5, forget plot] table [x=t, y=param_pe_ucb, col sep=comma] {figs/data_means2.csv};
            \addplot[name path=param_pe_lcb, line width=1pt, mark=none, red!30, opacity=0.5, forget plot] table [x=t, y=param_pe_lcb, col sep=comma] {figs/data_means2.csv};
            \addplot[red!30, forget plot, opacity=0.5] fill between [of=param_pe_lcb and param_pe_ucb];
            \addplot[line width=1pt, mark=none, red] table [x=t, y=param_pe_mean, col sep=comma] {figs/data_means2.csv};
            \addplot[black, thin, dashed] coordinates {(0, 0.99) (10, 0.99)};
            
            % Thin horizontal line at y = 0.8
            \addplot[black, thin, dashed] coordinates {(10, 0.8) (25, 0.8)};
            \addplot[black, thin, dashed] coordinates {(10, -17) (10, 2)};

            %PE constraint: (individual lines)
            % \foreach \i in {1,...,10}
            %     \addplot[line width=1pt, mark=none, red] table [x=t, y=param_pe_\i, col sep=comma] {figs/data_all.csv};

        \end{axis}

        %Proposed params:
        \begin{axis}[
            name = params_prop,
            at=(params_pe.below south), anchor=above north,
            width=\plotwidth, height=\plotheight,
            xlabel style={at={(axis description cs:0.5,-0.3)},anchor=north},
            ylabel={$A_{21}$},
            xlabel={Time step $k$},
            ylabel near ticks,
            xmin=0, xmax=25,
            ymin=-0.5, ymax=1.5,
            xtick={0,5,10,15,20,25},
            ytick={0,1},
            xmajorgrids=true,
            ymajorgrids=true,
            grid style=dashed,
            fill between/on layer=axis background,
            label style={font=\footnotesize, align=center},
            tick label style={font=\footnotesize} 
        ]      
            
            %Proposed: (confidence bounds)
            \addplot[name path=param_prop_ucb, line width=1pt, mark=none, green!30, opacity=0.5, forget plot] table [x=t, y=param_prop_ucb, col sep=comma] {figs/data_means2.csv};
            \addplot[name path=param_prop_lcb, line width=1pt, mark=none, green!30, opacity=0.5, forget plot] table [x=t, y=param_prop_lcb, col sep=comma] {figs/data_means2.csv};
            \addplot[green!30, forget plot, opacity=0.5] fill between [of=param_prop_lcb and param_prop_ucb];
            \addplot[line width=1pt, mark=none, green] table [x=t, y=param_prop_mean, col sep=comma] {figs/data_means2.csv};
            \addplot[black, thin, dashed] coordinates {(0, 0.99) (10, 0.99)};
            
            % Thin horizontal line at y = 0.8
            \addplot[black, thin, dashed] coordinates {(10, 0.8) (25, 0.8)};
            \addplot[black, thin, dashed] coordinates {(10, -17) (10, 2)};

            %Proposed: (individual lines)
            % \foreach \i in {1,...,10}
            %     \addplot[line width=1pt, mark=none, green] table [x=t, y=param_prop_\i, col sep=comma] {figs/data_all.csv};
            
        \end{axis}
        
    \end{tikzpicture}
    }
    %%%%%%%%%%%%%%%%%%%%%%%%%%%%%%%%%%%%%%%%%%%%%%%%%%%%%%%%%%%%%%%%%%%%%%%%%%%
    \subfigure[Volumes $V$.]{
    \begin{tikzpicture}
        % Nominal volumes:
        \begin{axis}[
            name = vol_nom,
            width=\plotwidth, height=\plotheighti,
            xlabel style={at={(axis description cs:0.5,-0.3)},anchor=north},
            ylabel={$V$},
            ylabel near ticks,
            xmin=0, xmax=25,
            ymin=-17, ymax=2,
            xtick={0,5,10,15,20,25},
            xticklabels={},
            ytick={-14,-7,-1},
            yticklabels={\tiny $10^{-14}$,  \tiny $10^{-7}$, \tiny $10^{-1}$},
            % ymode = log,
            xmajorgrids=true,
            ymajorgrids=true,
            grid style=dashed,
            fill between/on layer=axis background,
            label style={font=\footnotesize, align=center},
            tick label style={font=\footnotesize},
            title style = {font=\footnotesize}
        ]

            %Nominal: (confidence bounds)
            \addplot[name path=vol_nom_ucb, line width=1pt, mark=none, blue!30, opacity=0.5, forget plot] table [x=t, y=vol_nom_ucb, col sep=comma] {figs/data_means2.csv};
            \addplot[name path=vol_nom_lcb, line width=1pt, mark=none, blue!30, opacity=0.5, forget plot] table [x=t, y=vol_nom_lcb, col sep=comma] {figs/data_means2.csv};
            \addplot[blue!30, forget plot, opacity=0.5] fill between [of=vol_nom_lcb and vol_nom_ucb];
            \addplot[line width=1pt, mark=none, blue] table [x=t, y=vol_nom_mean, col sep=comma] {figs/data_means2.csv};            
            \addplot[black, thin, dashed] coordinates {(10, -17) (10, 2)};

            %Nominal: (individual lines)
            % \foreach \i in {1,...,10}
            %     \addplot[line width=1pt, mark=none, blue] table [x=t, y=vol_nom_\i, col sep=comma] {figs/data_all.csv};
            
        \end{axis}

        %PE volumes:
        \begin{axis}[
            name = vol_pe,
            at=(vol_nom.below south), anchor=above north,
            width=\plotwidth, height=\plotheighti,
            xlabel style={at={(axis description cs:0.5,-0.3)},anchor=north},
            ylabel={$V$},
            ylabel near ticks,
            xmin=0, xmax=25,
            ymin=-17, ymax=2,
            xtick={0,5,10,15,20,25},
            xticklabels={},
            ytick={-14,-7,-1},
            yticklabels={\tiny $10^{-14}$,  \tiny $10^{-7}$, \tiny $10^{-1}$},
            % ymode = log,
            xmajorgrids=true,
            ymajorgrids=true,
            grid style=dashed,
            fill between/on layer=axis background,
            label style={font=\footnotesize, align=center},
            tick label style={font=\footnotesize} 
        ]
            
            %PE constraint: (confidence bounds)
            \addplot[name path=vol_pe_ucb, line width=1pt, mark=none, red!30, opacity=0.5, forget plot] table [x=t, y=vol_pe_ucb, col sep=comma] {figs/data_means2.csv};
            \addplot[name path=vol_pe_lcb, line width=1pt, mark=none, red!30, opacity=0.5, forget plot] table [x=t, y=vol_pe_lcb, col sep=comma] {figs/data_means2.csv};
            \addplot[red!30, forget plot, opacity=0.5] fill between [of=vol_pe_lcb and vol_pe_ucb];
            \addplot[line width=1pt, mark=none, red] table [x=t, y=vol_pe_mean, col sep=comma] {figs/data_means2.csv};          
            \addplot[black, thin, dashed] coordinates {(10, -17) (10, 2)};

            %PE constraint: (individual lines)
            % \foreach \i in {1,...,10}
            %     \addplot[line width=1pt, mark=none, red] table [x=t, y=vol_pe_\i, col sep=comma] {figs/data_all.csv};

        \end{axis}

        %Proposed volumes:
        \begin{axis}[
            name = vol_prop,
            at=(vol_pe.below south), anchor=above north,
            width=\plotwidth, height=\plotheighti,
            xlabel style={at={(axis description cs:0.5,-0.3)},anchor=north},
            ylabel={$V$},
            xlabel={Time step $k$},
            ylabel near ticks,
            xmin=0, xmax=25,
            ymin=-17, ymax=2,
            xtick={0,5,10,15,20,25},
            ytick={-14,-7,-1},
            yticklabels={\tiny $10^{-14}$,  \tiny $10^{-7}$, \tiny $10^{-1}$},
            % ymode = log,
            xmajorgrids=true,
            ymajorgrids=true,
            grid style=dashed,
            fill between/on layer=axis background,
            label style={font=\footnotesize, align=center},
            tick label style={font=\footnotesize} 
        ]      
            
            %Proposed: (confidence bounds)
            \addplot[name path=vol_prop_ucb, line width=1pt, mark=none, green!30, opacity=0.5, forget plot] table [x=t, y=vol_prop_ucb, col sep=comma] {figs/data_means2.csv};
            \addplot[name path=vol_prop_lcb, line width=1pt, mark=none, green!30, opacity=0.5, forget plot] table [x=t, y=vol_prop_lcb, col sep=comma] {figs/data_means2.csv};
            \addplot[green!30, forget plot, opacity=0.5] fill between [of=vol_prop_lcb and vol_prop_ucb];
            \addplot[line width=1pt, mark=none, green] table [x=t, y=vol_prop_mean, col sep=comma] {figs/data_means2.csv};   
            \addplot[black, thin, dashed] coordinates {(10, -17) (10, 2)};

            %Proposed: (individual lines)
            % \foreach \i in {1,...,10}
            %     \addplot[line width=1pt, mark=none, green] table [x=t, y=vol_prop_\i, col sep=comma] {figs/data_all.csv};
            
        \end{axis}
        
    \end{tikzpicture}
    }
    \vspace{-2mm}
    \caption{Fault identification and control results for the passive (Pass.) and adaptive (Adapt.) RMPC, and our proposed (Prop.) \acr{AFD} strategy. In subfigure~a), the input trajectories of five exemplary runs are shown. In subfigures b) and~c), the solid lines represent the mean values over ten runs, and the shaded areas indicate the 95\% confidence intervals. The vertical dashed line marks the onset of Fault~2 at time step 10. 
}
    \label{fig:plots}
    \vspace{-5mm}
\end{figure*}

Fig.~\ref{fig:plots}b) shows the evolution of the parameter estimate $A_{12}$ as an example for all model parameters, which is 0.99 for $k \in \left[0,9\right]$ and drops to 0.8 thereafter, indicated by the dashed lines. The adaptive controller converges slowest to the true value due to a lack of excitation, while the passive controller could, in theory, estimate the parameter accurately as second fastest, but without an adaptation mechanism, this does not improve control performance. The proposed method achieves the fastest and most consistent convergence across trials, effectively adapting the model to the true system. Fig.~\ref{fig:plots}c) shows the corresponding uncertainty set volume over time. Similar trends emerge: the adaptive controller reduces the volume most slowly due to uninformative inputs, while the passive controller reduces it more quickly, though again without any control benefit. The proposed \acr{AFD} method shrinks the uncertainty set most rapidly following each fault, and once the volume drops below $V_\text{min}$, excitation ceases, allowing the controller to focus on tracking. 

Table~\ref{tab:my-table} summarizes key metrics for both fault scenarios. Best-performing results are highlighted in bold, detection times are computed for the cases where a fault is detected. For Fault~1, the proposed method achieves the fastest detection times, lowest final volumes, and highest detection rates, outperforming both the nominal and PE-constrained strategies. The passive controller performs second best and detects the more severe Fault~2 earlier due to immediate unanticipated system excitation. In contrast, the adaptive approach struggles most with FI; although increasing the PE gain could improve detection, it would degrade control performance. This highlights the need for a trade-off mechanism, such as the proposed $\beta$.

\begin{table}[t]
\caption{\footnotesize Identification results for Fault 1 (F1) and Fault 2 (F2).}
\vspace{-2mm}
\centering
\label{tab:my-table}
\footnotesize 
\begin{tabular}{ll|l|l|l|}
\cline{3-5}
                                                        &                        & \textbf{Passive} & \textbf{Adaptive} & \textbf{Proposed} \\ \hline
\multicolumn{1}{|c|}{\multirow{3}{*}{\textbf{F1}}} & Detection time         & 5.92       & 6.34       & \textbf{3.92 }      \\ \cline{2-5} 
\multicolumn{1}{|c|}{}                                  & Detection rate         & 84.2\%     & 44.2\%     & \textbf{99.2\%}    \\ \cline{2-5}
\multicolumn{1}{|c|}{}                                  & Final volume           & 8.9757e-10   & 1.9843e-04   &\textbf{3.8792e-10}   \\ \cline{2-5} 

\hline \hline
\multicolumn{1}{|l|}{\multirow{3}{*}{\textbf{F2}}} & Detection time         & \textbf{1.03}       & 1.75       & 1.33       \\ \cline{2-5} 
\multicolumn{1}{|l|}{}                                  & Detection rate         & \textbf{100\% }    & 96.7\%    & \textbf{100\% }   \\ \cline{2-5}
\multicolumn{1}{|l|}{}                                  & Final volume           &1.9458e-11  &8.5216e-08   &\textbf{1.4848e-13}  \\ \cline{2-5} 

\hline
\vspace{-9mm}
\end{tabular}
\end{table}

\section{Conclusion}
\label{sec:conclusion}
\vspace{-1mm}
Accurately identifying small, unknown faults remains a challenging task, particularly in systems that lack persistent excitation. While active identification methods improve sensitivity, handling unknown fault dynamics remains difficult. In this work, we proposed a robust and adaptive framework for \acr{AFD} that integrates FI with adaptive RMPC for a continuous set of unknown bounded faults. By penalizing the uncertainty volume via ellipsoidal set approximations and volume-dependent scaling, the method safely excites the system only when needed and gradually returns to nominal RMPC behavior as model confidence improves. This avoids permanent excitation and tuning sensitivity associated with PE-constrained methods. Simulation results show that the proposed strategy outperforms passive and adaptive PE-based RMPC approaches in detection speed, estimation accuracy, and closed-loop performance, while maintaining feasibility and satisfying constraints. It effectively balances control and exploration without requiring hard excitation constraints.

% A key enabler of fault detectability is the design of the initial uncertainty set, which must be sufficiently large to contain all relevant fault scenarios, yet stabilizable to ensure robust controller feasibility. 
In this context, finding an initial set for faults remains an open challenge. Further extensions will address reducing the computational complexity and incorporating scalable ellipsoidal approximations throughout the RMPC framework.

% \section*{ACKNOWLEDGMENT}

\vspace{-2mm}
\bibliography{mybib}
\bibliographystyle{ieeetr}

\end{document}